\documentstyle[12pt]{article}
\begin{document}
\title{   Convergence of Fourier series at or  beyond endpoint
}
\author{{   Shunchao Long }\\
 }
\date{}
\maketitle
\begin{center}
\begin{minipage}{120mm}
\vskip 0.1in
{
\begin{center}
{\bf Abstract}
\end{center}
~~  We consider several problems at or beyond endpoint in harmonic
analysis. The solutions of these problems  are related to the
estimates of some classes of sublinear operators. To do this, we
introduce  some new functions spaces $RL^{p,s }_{|x|^{\alpha}}({\bf
R}^n)$ and $\dot{R}L^{p,s }_{|x|^{\alpha}}({\bf R}^n)$, which play
an analogue role  with the classical Hardy spaces $H^p({\bf R}^n)$.
These spaces are subspaces of $L^p_{|x|^{\alpha}}({\bf R}^n)$ with
$1<s<\infty, 0<p\leq s$ and $-n<\alpha<n(p-1)$, and $\dot{R}L^{p,s
}_{|x|^{\alpha}}({\bf R}^n)  \supset L^s({\bf R}^n)$ when $
-n<\alpha<n(p/s-1)$. We prove the following results.

\par  ~~ First,    $\mu_\alpha$-a.e. convergence and ${L}^{p}_{|x|^{\alpha}}({\bf R} )$ -norm convergence  of Fourier
 series hold for all functions in  $ RL^{p,s }_{|x|^{\alpha}}({\bf R} )$  and  $ \dot{R}L^{p,s
}_{|x|^{\alpha}}({\bf R} )$ with  $1<s<\infty, 0<p\leq s$ and
$-1<\alpha<p-1$, where $\mu_\alpha(x)=|x|^{\alpha}$;

 \par ~~ Second,   many  sublinear operators
 initially defined for the functions  in $L^p({\bf R}^n)$ with   $1<p<\infty$, such as Calder\'{o}n-Zygmund
operators,  C.Fefferman's singular multiplier operator,
R.Fefferman's singular integral operator, the Bochner-Riesz means at
the critical index,  certain oscillatory singular integral
operators, and so on, admit extensions which map
 $RL^{p,s }_{|x|^{\alpha}}({\bf R}^n)$ and $\dot{R}L^{p,s
}_{|x|^{\alpha}}({\bf R}^n)$ into $L^p_{|x|^{\alpha}}({\bf R}^n)$
with $1<s<\infty, 0<p\leq s$ and $-n<\alpha<n(p-1)$;
\par ~~ Final,     Hardy-Littlewood maximal
 operator is  bounded  from $RL^{p,s }_{|x|^{\alpha}}({\bf R}^n)$ (or $\dot{R}L^{p,s
}_{|x|^{\alpha}}({\bf R}^n)$) to ${L}^{p}_{|x|^{\alpha}}({\bf R}^n)$
for $ 1<s<\infty$ and $0<p\leq s$  if and only if
$-n<\alpha<n(p-1)$.

}
\end{minipage}
\end{center}
\vskip 0.1in
\baselineskip 16pt

\par \begin{center}
{\bf Table of contents}
\end{center}
\par \textrm{1.}
 Introduction
\par \textrm{2.}
  $RL^{p,s}_{|x|^\alpha}$ and $
\dot{R}L^{p,s}_{|x|^\alpha} $: definitions  and basic properties
\par \textrm{3.}
  Some classes of  operators
\par \textrm{4.}
     Hardy-Littlewood maximal operator
\par \textrm{5.}  Singular integrals
\par \textrm{6.} Fourier series

\par \section*
 {\bf  1 Introduction  }

\par~~~~In this paper we study  several fundamental  problems at or beyond endpoint in harmonic
analysis, such as  convergence of Fourier series on ${\bf R} $,
extension of some singular integrals and boundedness of
Hardy-Littlewood maximal operator on ${\bf R}^n$. These problems are
known to be soluble on $L^p({\bf R}^n)$, ($n=1$ or $n\geq 1$), for
$1<p<\infty$, but not for $0<p\leq 1$, even not on $H^p({\bf R}^n)$.
 One naturally asks how these known results for
 $1<p<\infty$ are extend to $0<p\leq 1$.
 This can be done when we consider the
 problems on the weighted Lebesgue spaces $L^p_{|x|^\alpha}({\bf R}^n)$.
 We show  that these problems are soluble on some
 subspaces, $ RL^{p,s}_{|x|^\alpha}({\bf R}^n)$
 and
 $ \dot{R}L^{p,s}_{|x|^\alpha}({\bf R}^n)$, of $L^p_{|x|^\alpha}({\bf R}^n)$ with  $0<p\leq 1,-n
 <\alpha<n(p-1)$ and $1<s<\infty$.
These spaces  are extensions of the classical Hardy
 spaces and the block spaces of M.Taibleson and G.Weiss. They play an analogue role  with the
classical Hardy spaces  in this paper. And $
\dot{R}L^{p,s}_{|x|^\alpha}({\bf R}^n)$ includes $L^s ({\bf R}^n) $
with $1<s<\infty, 0<p\leq s $ and $ -1<\alpha< p/s-1.$

   \par We denote by
$S_Nf$   the partial sums of the Fourier Series of a function $f$ on
$T=[-\pi,\pi]$ or the Fourier integral of $f$ on ${\bf R}$ .
   A famous result of L.Carleson \cite{Carl} states that:  if $f\in L^2(T)$ then
$S_Nf(x)$ converges to $ f(x)$ almost everywhere as $N\rightarrow
\infty$. R.Hunt \cite{Hunt} extended this result to $1<p<\infty$.
Alternative proofs of L.Carleson's theorem were provided by
C.Fefferman \cite{F3} and M.Lacey and C.Thiele  \cite{LT}. In fact,
M.Lacey and C.Thiele proved the theorem on the line {\bf R}. An
alternative proof of Hunt's theorem was provided by
 L.Grafakos, T.Tao and
E.Terwilleger  \cite{GTT}, and they proved the theorem on the line
{\bf R}.

\par The convergence problem of Fourier series for the endpoint $p=1$ was
considered by many mathematicians.  In 1923, A.Kolmogorov \cite{Ko1}
proved that there exists a $f\in L^1$ whose Fourier series diverges
at almost every point. This shows that Carleson-Hunt's theorem does
not hold when $p=1$.
  In fact, the function $f$ constructed by  Kolmogorov is in $
L \log \log L(T) \subset L^1(T)$.
 Kolmogorov's result was extended to $f\in L (\log \log L)^{1-\varepsilon}(T)
 $ by Y.M.Chen \cite{Chen}, and to $f\in L \varphi(L)(T)$
  where $\varphi$ satisfies $\varphi(t)=o((\log t/\log \log
  t)^{1/2})$
  by S.V.Konyagin \cite{Kon2}.
   The example  of Kolmogorov can also be modified to yield a function $f\in H^1(T)$ whose Fourier
series diverges at almost every point (see   \cite{Zy}, Chapter 8).

\par Some positive results of extension of Carleson-Hunt theorem to $p=1$ have been
obtained by several authors. R.Hunt \cite{Hunt} proved   that if
$f\in L{\rm log} L{\rm  log  } L (T)$   then $S_Nf(x)$ converges to
$f(x)$ almost everywhere. This result was sharpened to $f\in L{\rm
log} L{\rm  log log } L (T)$ by P.Sjolin \cite{Sj},  and to $f\in
L{\rm log} L{\rm log log log } L (T)$ by N.Y.Antonov \cite{An}. Also
there are some quasi-Banach spaces of functions with   a.e.
convergent Fourier series, see \cite {S3,S1} and \cite{Ar2}. And
M.Taibleson and G.Weiss \cite{TW1} proved that if $f\in B_q$
(generated by $q$-blocks) $\subset L^1(T),$ then $S_Nf(x)$ converges
to $ f(x)$ almost everywhere.

\par   Our first result of this paper is to extend  Carleson-Hunt's
theorem  to $0<p\leq 1$ as follows: let  $1< s<\infty, 0<p\leq s$
 and $
-1<\alpha< p-1.$ Then
 $$\lim_{N\rightarrow \infty }S_Nf(x)= f(x),~~~~\mu_\alpha-{\rm a.e.}$$
 for all  $f\in   RL^{p,s}_{|x|^\alpha}({\bf R})$
 and
 $ \dot{R}L^{p,s}_{|x|^\alpha}({\bf R})$, where
 $\mu_\alpha(x)=|x|^\alpha$.

\par On the other hand, an early celebrated result due to M.Riesz \cite{Riesz} states that
$S_Nf(x)$ converges to $ f(x)$   in $L^p({\bf R})$-norm   for
$1<p<\infty. $ But for $p=1,$ Kolmogorov's example shows that there
exists a $f\in L^1({\bf R})$ whose Fourier series does not converge
to $f$ in $L^1({\bf R})$. And the Fourier series of the function
$f\in H^1({\bf R})$ mentioned above, yielded by the example of
Kolmogorov,  does  not converge to $f$ in $L^1({\bf R})$ as well.

\par Here, our second result is to extend the norm convergence of Fourier series  to $0<p\leq 1$ as
follows: let  $1< s<\infty, 0<p\leq s$
 and $
-1<\alpha< p-1,$ then
$$\lim_{N\rightarrow \infty } \|S_Nf- f\|_{L^p_{|x|^{\alpha}}({\bf R})}=0$$
 for all  $f\in   RL^{p,s}_{|x|^\alpha}({\bf R})$
 and
 $ \dot{R}L^{p,s}_{|x|^\alpha}({\bf R})$.

\par The partial sums operator $S_N$ of the Fourier Series or the Fourier integral
is  a sublinear operator which is closely related to Hilbert
transform $H$ since the identity (6.2) below.
The convergence results above rely on some estimates of $S_N$. In
fact, in this paper we shall consider some classes of sublinear
operators which are related to $S_N$, many singular integral
operators and Hardy-Littlewood maximal operator.
\par It is well-known that the theory for singular integral operators
plays an important role in harmonic analysis. Many classical
singular integral operators, such as Calder\'{o}n-Zygmund singular
integral operators,   C.Fefferman's strong singular multiplier
operator, R.Fefferman's singular integral operator,  some
oscillatory singular integral operator, and so on, are initially
defined on Schwartz functions space. These operators are known to
extend boundedness linear operators on $L^p({\bf R}^n)$ for
$1<p<\infty$, but not for $0<p\leq 1$. It is  known that
Calder\'{o}n-Zygmund singular integral operators admit the
extensions which map $H^p({\bf R}^n)$ into $L^p({\bf R}^n)$ for
$p\leq 1$. However, some of these operators, such as certain
C.Fefferman's strong singular multiplier operators and certain
oscillatory singular integral operators, fail to  map $H^p({\bf
R}^n)$ into $L^p({\bf R}^n)$, or $H^p({\bf R}^n)$ into itself, for
certain $0<p\leq 1$. (See \cite{St1970, FS72, Sjolin76, Pan3,
Pan6}).


\par In this paper, our third work is to show that these singular integral
operators $T$ mentioned above admit some extensions which map $
RL^{p,s}_{|x|^\alpha}({\bf R}^n)$
 (or
 $ \dot{R}L^{p,s}_{|x|^\alpha}({\bf R}^n)$) into $L^p_{|x|^{\alpha}}({\bf R}^n)
 $ with $1<s<\infty, 0<p\leq
s$ and $-n<\alpha<n(p-1)$, and satisfy
$$Tf=\sum_{j=1}^\infty \lambda_jTa_j, ~~\mu_\alpha-{\rm a.e.}, $$
  for all $f=\sum_{j=1}^\infty
\lambda_ja_j  $ where each $a_i$ is a central R-$( p,s,
\alpha)-$block and $\sum_{j=1}^\infty |\lambda_j|^{\bar{p}}<\infty $
with  $  \bar{p} = {\rm min} \{p, 1\}  .$

It is also known that the Hardy-Littlewood maximal operator $M$ is
bounded on  $L^p({\bf R }^n)$ for $1<p<\infty$. Due to this property
and its character  of controlling many operators  $M$ is very useful
in Harmonic analysis. So, it is also important to extend its
$L^p({\bf R }^n)$-boundedness property for $1<p<\infty$  to $0<p\leq
1$.
 Unfortunately,  $M $ is not a bounded operator
on  $L^1({\bf R }^n)$, even not   from $H^p({\bf R }^n)$ into
$L^p({\bf R }^n)$ for $0<p\leq 1$, although there is   a substitute
result, namely that   it is bounded from $L^1({\bf R }^n)$ into
weak-$L^1({\bf R }^n)$. In fact, there is no  non-zero subspace in
$L^p ({\bf R }^n)$ from which
 $Mf$ is  bounded
 since it is never in  $L^p ({\bf R }^n)$ unless $f=0$ for $0<p\leq 1$.
(See \cite{Grafakos}).

\par Surprisingly, we find here that the weighted norm inequality of
$M$ can be extended to $0<p\leq 1.$ Our fourth  result is that $M$
maps $ RL^{p,s}_{|x|^\alpha}({\bf R}^n)$
 (or
 $ \dot{R}L^{p,s}_{|x|^\alpha}({\bf R}^n)$) into $L^p_{|x|^{\alpha}}({\bf R}^n)
 $ with $1<s<\infty, 0<p\leq
s$ and $-n<\alpha<n(p-1)$.

\par In fact, we have more for Hilbert transform and Hardy-Littlewood
maximal operator.
 According to the famous results by B.Muckenhoupt \cite{Muck} and
R.Hunt, B.Muckenhoupt and R.Wheeden \cite{HMW}, it is known that,
for $1<p<\infty$,
$$
\|Tf\|_{L^{p}_{|x|^{\alpha}}({\bf R }^n)} \leq C
\|f\|_{L^{p}_{|x|^{\alpha}}({\bf R }^n)}
$$
if and only if $-n  <\alpha<n(p-1 ),  $ where $T$ are Hilbert
transform and Hardy-Littlewood maximal operator.

\par Here,  we extend these to $0<p\leq 1$ as follows: let $1<
s<\infty,0<p \leq s,$ then
$$ \|Tf\|_{L^{p}_{|x|^{\alpha}}({\bf R }^n)} \leq C
\|f\|_{RL^{p,s}_{|x|^{\alpha}}({\bf R }^n)} {\rm ~~or ~~}
\|Tf\|_{L^{p}_{|x|^{\alpha}}({\bf R }^n)} \leq C
\|f\|_{\dot{R}L^{p,s}_{|x|^{\alpha}}({\bf R }^n)}
$$
if and only if $-n  <\alpha<n(p-1 ). $

\par The paper is organized as follows. In section 2, we introduce
the spaces $RL^{p,s}_{|x|^\alpha}({\bf R }^n) $ and $
\dot{R}L^{p,s}_{|x|^\alpha}({\bf R }^n) $, and  list some of their
properties. In this section, we see that the spaces $
\dot{R}L^{p,s}_{|x|^\alpha}({\bf R }^n) $ are larger than $L^s$ for
all $1<s<\infty$. In section 3 we discuss  the boundedness
properties of some classes of sublinear
 operators.
 The boundedness of Hardy-Littlewood maximal operator is discussed   in section
 4.   The extension properties of some singular integral operators are obtained  in
 section 5. Section 6 is devote to   the convergence of Fourier series,
 including pointwise convergence and norm convergence.

\section*
 {\bf  2   $RL^{p,s}_{|x|^\alpha}$ and $
\dot{R}L^{p,s}_{|x|^\alpha} $: definitions  and basic properties
}

\par The spaces $
BL^{p,s}_{|x|^\alpha} $ introduced in \cite{Longs2,Longs}  are very
useful in the study for norm convergence of Bochner-Riesz means and
the spherical means at or beyond endpoint. Unfortunately, we can not
obtain new pointwise convergence result of Fourier series by using $
BL^{p,s}_{|x|^\alpha} $. Because  a restriction $n(p/s-1) <\alpha <
n(p-1)$  is needed in our use of $ BL^{p,s}_{|x|^\alpha} $ in
\cite{Longs} , and $ BL^{p,s}_{|x|^{\alpha}}\subset
L^{\frac{np}{n+\alpha}}$ at this time. Here, we introduce some new
functions spaces $RL^{p,s}_{|x|^\alpha}$ and $
\dot{R}L^{p,s}_{|x|^\alpha} $ which are slightly different from the
spaces $B{L}^{p,s }_{|x|^{\alpha}}$. These spaces generate by the
annulus blocks centered at origin. So, they are also extensions of
the classical Hardy spaces and Taibleson and Weiss's blocks spaces.
\par
 Denote
 $B_k=\{x\in {\bf R}^n: |x| \leq 2^k \},  C_k=B_k\setminus B_{k-1}
 $ for $k\in {\bf Z},$ and $ \widetilde{C}_k= C_k$ for $k\in {\bf N}$ and $ \widetilde{C}_0=
 B_0$.
  $\chi
_{E}$ is the characteristic function of set $E$.

\par {\bf Definition 1} ~~Let $
0< s \leq \infty, 0< p < \infty,- \infty < \alpha < \infty .$
\par A. A function $a(x)$ is said to be a central R-$ (p, s, \alpha)$-block on ${\bf R}^n$, if
\par (i)~~~~supp $a\subseteq C_k\subset {\bf R}^n,k\in {\bf Z},$
\par (ii)~~~~$\|a\|_{L^{s}   }\leq |B_k)|^{-\alpha/pn-1/p +1/s};$

\par B. A function $a(x)$ is said to be a central R-$ (p, s, \alpha)$-block of restrict type on ${\bf R}^n$, if
\par (i)~~~~supp $a\subseteq \widetilde{C}_k\subset {\bf R}^n,k\in {\bf N}\cup\{0\},$
\par (ii)~~~~$\|a\|_{L^{s}   }\leq |B_k|^{-\alpha/pn-1/p +1/s}$.

\par {\bf Definition 2}~~ Let $0< s \leq \infty, 0<p< \infty ,
 - \infty < \alpha < \infty .$
\par A. The homogeneous functions spaces
$ RL^{p,s}_{|x|^{\alpha}}({\bf R}^n)$ is defined as
\begin{eqnarray*}
 RL^{p,s}_{|x|^{\alpha}}({\bf R}^n)=\{  f&: &
 f=\sum\limits_{k=-\infty}^{\infty} \lambda _ka_k ,
 \\&&\textrm { where
 each $a_k$ is a central R-$(p,s, \alpha)$-block on ${\bf R}^n$,}
 \\&& \sum\limits_{k=-\infty}^{\infty} |\lambda _k|^{\bar{p}} <+ \infty \},
 \end{eqnarray*}
  the
 "convergence" is meant in the sense of $\mu_\alpha$-a.e. convergence.
 Moreover, we define a quasinorm on $ RL^{p,s}_{|x|^{\alpha}}({\bf
 R}^n)$ by
$$\|f\|_{RL^{p,s}_{|x|^{\alpha}}({\bf R}^n)}=
 \inf \left(\sum\limits_{k=-\infty}^{\infty}|\lambda _k|^{\bar{p}}\right)^{1/{\bar{p}}},$$
where the infimum is taken over all the decompositions of $f$ as
above.
\par B. The non-homogeneous functions spaces $\dot{R}L^{p,s}_{|x|^{\alpha}}({\bf R}^n)$ are defined  replacing the
central R-$(p,s, \alpha)$-blocks in the definition above by central
R-$ (p, s, \alpha)$-blocks of restrict type.

\par Simply, we denote $RL^{p,s}_{|x|^{\alpha}}({\bf R}^n)$ and $\dot{R}L^{p,s}_{|x|^{\alpha}}({\bf R}^n)$ by
${R}L^{p,s}_{|x|^{\alpha}}$ and $\dot{R}L^{p,s}_{|x|^{\alpha}}$
respectively as $n>1$.

\par Now, let us state a number of properties of $
RL^{p,s}_{|x|^{\alpha}} $ and $ \dot{R}L^{p,s}_{|x|^{\alpha}}. $

  \par {\bf Proposition 2.1}  Let $1< s< \infty, 0<p\leq s,
-n<\alpha <n(p-1)$.  Then
$$
RL^{p,s}_{|x|^{\alpha}}  \subset {L}^{p}_{|x|^{\alpha}} {\rm~ and ~}
\dot{R}L^{p,s}_{|x|^{\alpha}}  \subset {L}^{p}_{|x|^{\alpha}}.
$$

  \par {\bf Proposition 2.2} Let $0<p\leq s\leq\infty$. Then,
\par i. if $ n(p/s-1) \leq \alpha<\infty$,
$$RL^{p,s}_{|x|^{\alpha}}\subset L^{\frac{np}{n+\alpha}}
 {\rm  ~ and ~}
\dot{R}L^{p,s}_{|x|^{\alpha}}\subset L^{\frac{np}{n+\alpha}},$$

 \par ii. if $-n<\alpha \leq n(p/s-1)$,
$$ L^s\subset \dot{R}L^{p,s}_{|x|^{\alpha}},$$

   \par iii. if $ \alpha = n(p/s-1)$,
$$ L^s = \dot{R}L^{p,s}_{|x|^{\alpha}}. $$

 \par {\bf Proposition 2.3} Let  $1<s\leq \infty,0<p\leq s,-n <\alpha< n(p-1)$. Then, $
RL^{p,s}_{|x|^{\alpha}}, $ and $ \dot{R}L^{p,s}_{|x|^{\alpha}} $ are
complete normed-spaces.

 \par {\bf Proposition 2.4}
Let $1\leq s<\infty, 0<p<\infty,-n<\alpha<\infty. $ Then
$C_0^\infty$ and $\mathcal{S} $ are dense in  $
R\dot{L}^{p,s}_{|x|^{\alpha}} $ and $ R{L}^{p,s}_{|x|^{\alpha}} $.

\par Proposition 2.1 is a corollary of Theorem 4.1 below.
 See \cite{Longs2} for  proofs of Propositions 2.2-2.4.

\par \section*
 {\bf  3 Some classes of  operators  }
 \par
Let  operators $T_1, T_2$ and $ T_0 $ be defined for any integrable
function with compact support and  satisfy the size conditions
below, respectively,
$$|T_1f(x)| \leq C \|f\|_{L^{1}}/|x|^n,\eqno {(3.1)}$$
when supp $f\subseteq C_k$ and $|x| \geq 2^{k+1}$ with $k\in {\bf
Z};$
$$|T_2f(x)| \leq C 2^{-kn}\|f\|_{L^{1}},\eqno {(3.2)}$$
when supp $f\subseteq C_k$ and $|x| \leq 2^{k-2}$ with $k\in {\bf
Z};$ and
$$T_0 ~{\rm satisfies~ (3.1)~
 and ~(3.2).}$$

{\bf Theorem 3.1} Let $1\leq s\leq \infty,0<p\leq p_0\leq s ,
-n<\alpha <n(p-1)$.
  Suppose that   $T_i, i=0,1,2$, defined as above,  are bounded on $
{L^{p _0} }$.
  For all central R-$(  p,s,
\alpha)-$blocks  $a$, we have,

\par A) if   $-n  < \alpha < n(p/p _0-1  ) ,$ then
$$
\|T_2a\|_{L^{p}_{|x|^{\alpha}}}\leq C ;$$
\par B) if  $ n(p/p _0-1  )< \alpha < n(p-1  ) ,$
then
$$ \|T_1a\|_{L^{p}_{|x|^{\alpha}}}\leq C
;$$
\par C) if   $- n  < \alpha < n(p-1 ),$ then
$$ \|T_0a\|_{L^{p}_{|x|^{\alpha}}}\leq C
,$$
 where  $C$ are independent of
$a$.

\par {\bf Proof}  Let $T$ be the  operators above, $a$ be a central R-
$(p,s,\alpha)$-block with  supp$a\subseteq C_{k_0}$ and
$\|a\|_{L^s}\leq \| B_{k_0}
\|^{-\frac{\alpha}{np}-\frac{1}{p}+\frac{1}{s}}$. For $ 0<p<\infty,$
we see that
\begin{eqnarray*}
\|Ta\|^p_{L^{p}_{|x|^{\alpha}}}
\approx \sum\limits_{j=-\infty}^{\infty}|B_j|^{\alpha /n}\|(Ta)\chi
_{C_j}\|^p_{L^{p}}.
\end{eqnarray*}
 Let
$$
I_{1}=\sum\limits_{j=-\infty}^{k_0-2}|B_j|^{\alpha  /n}\|(Ta)\chi
_{C_j}  \|^p_{L^{p}},
$$
$$
I_{3}=\sum\limits_{j=k_0+2 }^{\infty}|B_j|^{\alpha  /n}\|(Ta)\chi
_{C_j}  \|^p_{L^{p}},
$$
$$
I_{2A}=\sum\limits_{j=k_0-1 }^{\infty}|B_j|^{\alpha  /n}\|(Ta)\chi
_{C_j}  \|^p_{L^{p}},
$$
$$
I_{2B}=\sum\limits_{j=-\infty}^{k_0+1}|B_j|^{\alpha /n}\|(Ta)\chi
_{C_j}  \|^p_{L^{p}},
$$
$$
I_{2C}=\sum\limits_{j=k_0-1}^{k_0+1}|B_j|^{\alpha  /n}\|(Ta)\chi
_{C_j}  \|^p_{L^{p}}.
$$
It is easy to see that
$$
 \|Ta\|^p_{L^{p}_{|x|^{\alpha}}} \approx I_{1}+I_{2A};
 ~~\|Ta\|^p_{L^{p}_{|x|^{\alpha}}} \approx I_{2B}+I_{3};~~
  \|Ta\|^p_{L^{p}_{|x|^{\alpha}}} \approx I_{1}+I_{2C}+I_3.
$$
For $I_1,$ we use the facts that $k_0\geq j+2$ and $x\in C_j.$ By
(3.2) and H\"{o}lder's inequality,  we have
$$
|T_2a(x)|\leq C2^{-k_0n}
 \|a\|_{L^{1}}\leq C2^{-k_0n}\|a\|_{L^s }|B_{k_0} |^{ 1/s'}
\leq C |B_{k_0} |^{ -\alpha  /np-1/p} ,
$$
where $s'$ such that $\frac{1}{s}+\frac{1}{s'}=1$ for $1\leq s \leq
\infty$
It follows that
\begin{eqnarray*}
 I_{1}&\leq &
C\sum\limits_{j=-\infty}^{k_0-2}|B_j|^{\alpha /n} |B_{k_0} |^{
-\alpha  /n-1 } \|\chi _{C_j}  \|^p_{L^{p}}
   \leq  C\sum\limits_{j=-\infty}^{k_0-2}
\left( \frac{|B_j|}{|B_{k_0}|}\right)^{\alpha  /n+1}=C
\end{eqnarray*}
when $ \alpha/n+1 >0$. For $I_3,$ we use the facts that $k_0\leq
j-2$ and $x\in C_j .$ By (3.1), we have
\begin{eqnarray*}
|T_1a(x)|\chi _{C_j}  (x) &\leq& C2^{-jn}
 \|a\|_{L^{1}}
\leq C2^{-jn}
 \|a\|_{L^s }|B_{k_0}|^{ 1/s'}
\\
&\leq& C 2^{-jn} |B_{k_0} |^{ 1-\alpha  /np-1/p} .
\end{eqnarray*}
It follows that
\begin{eqnarray*}
 I_{3}&\leq &
C\sum\limits_{j=k_0+2 }^{\infty}|B_j|^{\alpha  /n} 2^{-jnp}
 |B_{k_0}|^{ p-\alpha  /n-1}
\|\chi _{C_j}  \|^p_{L^{p}}
 \leq
C\sum\limits_{j=k_0+2 }^{\infty} \left(
\frac{|B_j|}{|B_{k_0}|}\right)^{\alpha  /n+1-p}=C
\end{eqnarray*}
when $ \alpha+n-np<0$.

For $I_{2A}, I_{2B}$¡¡and $ I_{2C},$ let us first  estimate
$|B_j|^{\alpha /pn}\|\left( Ta  \right)\chi _{C_j}\|_{L^{p}}$.
 Let $0<p<\infty, p\leq p_0\leq s\leq \infty$. Using H\"{o}lder's inequality twice
 and $L^{p_0}$ - boundedness of $T$, we have:
\begin{eqnarray*}
|B_j|^{\alpha /pn}\|\left( Ta  \right)\chi _{C_j}\|_{L^{p}} & \leq &
|B_j|^{\alpha /pn} \left( \|\left( Ta   \right)\chi _{C_j}
 \|^p_{L^{p_0} } \|\chi _{C_j}   \|_{L^{(p_0/p)'} } \right)^{1/p}
\\
 & \leq & C
|B_j|^{\alpha /pn}
 \| a
\| _{L^{p_0} } |C_j|^{1/p-/p_0}
\\
 & \leq & C
|B_j|^{\alpha /pn}
 \| a
\| _{L^{s} }|B_{k_0}|^{ 1/p_0-1/s} |C_j|^{1/p-/p_0}
\\
& \leq & C \left(\frac{ |B_j|}{ |B_{k_0}|}\right)^{\alpha
/pn+1/p-1/p_0}
 .
\end{eqnarray*}
It follows that
\begin{eqnarray*}
I_{2A} &\leq & C\sum\limits_{j=k_0-1 }^{\infty}
\left(\frac{|B_j|}{|B_{k_0}|} \right)^{(\alpha /pn+1/p-1/p_0)p}
 =C
 \end{eqnarray*}
when $ \alpha<pn(1/p_0-1/p)$;
\begin{eqnarray*}
I_{2B} &\leq & C \sum\limits _{j=-\infty}^{k_0+1} \left(\frac{
|B_j|}{ |B_{k_0}|}\right)^{(\alpha /pn+1/p-1/p_0)p}=C
 \end{eqnarray*}
 when $ \alpha>pn(1/p_0-1/p)$;
and
\begin{eqnarray*}
I_{2C} &\leq & C \sum\limits _{j=k_0-1}^{k_0+1} \left(\frac{ |B_j|}{
|B_{k_0}|}\right)^{(\alpha /pn+1/p-1/p_0)p}=C
 \end{eqnarray*}
when $-\infty <\alpha<+\infty.$ Then, the desired results  can be
obtained  from the estimates above, in fact,  B) from $I_3$ and
$I_{2B}$, A) from $I_1$ and $I_{2A}$, C) from $I_1, I_3$ and
$I_{2C}$. Thus, we finish the proof of Theorem 3.1

\par From Theorem 3.1, we have

\par
{\bf Theorem 3.2} Let $1\leq s\leq \infty,0<p\leq p_0\leq s ,
-n<\alpha <n(p-1)$. Suppose  that  $T_i, i=0,1,2$, defined as above,
are bounded on $ {L^{p _0} }$,
  and
$$
|T_if|\leq \sum|\lambda_j||T_ia_j| , ~~\mu_\alpha-{\rm a.e.},
\eqno(3.3)
$$
for  $f=\sum\lambda_ja_j \in RL^{p,s} _{|x|^{\alpha}}$, where each
$a_j$ is a central R-$(  p,s, \alpha)-$block. We have,

\par A) if   $-n  < \alpha < n(p/p _0-1  ) ,$ then
$$
\|T_2f\|_{L^{p}_{|x|^{\alpha}}}\leq C
\|f\|_{RL^{p,s}_{|x|^{\alpha}}};$$
\par B) if  $ n(p/p _0-1  )< \alpha < n(p-1  ) ,$
then
$$ \|T_1f\|_{L^{p}_{|x|^{\alpha}}}\leq C
\|f\|_{RL^{p,s}_{|x|^{\alpha}}};$$
\par C) if   $- n  < \alpha < n(p-1 ),$ then
$$ \|T_0f\|_{L^{p}_{|x|^{\alpha}}}\leq C
\|f\|_{RL^{p,s}_{|x|^{\alpha}}},$$
 where  $C$ are independent of
$f$.

\par {\bf Proof}
Let $f\in RL^{p,s} _{|x|^{\alpha}}$, then for any $\varepsilon
>0,$ there exists $f=\sum\lambda_ja_j$ where each $a_j$ is a central
R-$(  p,s, \alpha)-$block such that
 $$\left( \sum
|\lambda_j|^{\bar{p}}\right)^{1/\bar{p}}\leq \|f\|_{RL^{p,s}
 _{|x|^{\alpha}}}+\varepsilon. $$
Then, by (3.3), Minkowski inequality and Theorem 3.1,   we have
 $$\|T_i f\|^p_{L^{p} _{|x|^{\alpha}}}\leq
\sum |\lambda_j|^p\|T_i a_j\|^p_{L^{p} _{|x|^{\alpha}}} \leq C \sum
|\lambda_j|^p \leq C \|f\|^p_{RL^{p,s}
 _{|x|^{\alpha}}}+\varepsilon^p
 $$
when $0<p\leq 1 $, and
 $$\|T_i f\|_{L^{p} _{|x|^{\alpha}}}\leq
\sum |\lambda_j|\|T_i a_j\|_{L^{p} _{|x|^{\alpha}}} \leq C \sum
|\lambda_j| = C \left( \sum |\lambda_j|^{\bar{p}}\right)^{1/\bar{p}}
\leq C\|f\|_{RL^{p,s}
 _{|x|^{\alpha}}}+\varepsilon
 $$
when $p>1 $. Since $\varepsilon$ is arbitrary, we have $\|T_i
f\|_{L^{p} _{|x|^{\alpha}}} \leq C\|f\|_{RL^{p,s}
 _{|x|^{\alpha}}} .
 $ Thus, we finish the proof of Theorem 3.2.

\par For the   non-homogeneous spaces $\dot{R}L^{p,s}_{|x|^{\alpha}}$
, we have a similar theorem whose proof is similar to that of
Theorems 3.1 and 3.2.

 \par Let $\overline{T}_1, \overline{T}_2, \overline{T}_0 $ be sublinear
operators saitisfing the size conditions below respectively,
$$| \overline{T}_1f(x)| \leq C \|f\|_{L^{1}}/|x-x_0|^n,\eqno {(3.4)}$$
when supp $f\subseteq B_0(x_0)$ and $|x-x_0| > 2$ or when supp
$f\subseteq C_k(x_0)$ and $|x-x_0| \geq 2^{k+1}$ with $k\in  {\bf
N},$ and
$$|\overline{T}_2f(x)| \leq C 2^{-kn}\|f\|_{L^{1}},\eqno {(3.5)}$$
when supp $f\subseteq C_k(x_0)$ and $|x-x_0| \leq 2^{k-2}$ with
$k\in {\bf N}\setminus\{1\},$ and
$$\overline{T}_0 ~{\rm satisfies~ (3.4)~
 and~ (3.5).}$$

{\bf Theorem 3.3}~~ Let $1 \leq s \leq \infty, 0<p\leq p_0\leq s ,
-n<\alpha <n(p-1)$.
  Suppose  that  $\overline{T}_i, i=0,1,2$, defined as above,  are bounded on $
{L^{p _0} }$,   and  satisfy (3.3) for  $f=\sum\lambda_ja_j \in
\dot{R}L^{p,s} _{|x|^{\alpha}},$  we have
\par A) if   $-n  < \alpha < n(p/p _0-1  ) ,$ then
 $$
\|\overline{T}_1f\|_{L^{p}_{|x|^{\alpha}}}\leq C
\|f\|_{\dot{R}L^{p,s}_{|x|^{\alpha}}};
$$
\par B) if  $ n(p/p _0-1  )< \alpha < n(p-1  ) ,$  then
$$
\|\overline{T}_2f\|_{L^{p}_{|x|^{\alpha}}}\leq C
\|f\|_{\dot{R}L^{p,s}_{|x|^{\alpha}}};
$$
\par C) if   $- n  < \alpha < n(p-1 ),$ then
 $$ \|\overline{T}_0f\|_{L^{p}_{|x|^{\alpha}}}\leq C
\|f\|_{\dot{R}L^{p,s}_{|x|^{\alpha}}},
$$
 where $C$ are independent of $f$.

{\bf Corollary 3.1} Let $1\leq s < \infty,0<p\leq s,- n < \alpha <
n(p-1)  $.
 Suppose  that a  operator $T$ satisfies the size condition
$$|Tf(x)|\leq C \int_{{\bf R}^n}\frac{|f(y)|}{|x-y|^n}dy, ~~~~ x\notin
{\rm supp}f, \eqno{(3.6)}$$
for any integrable function with compact support and  $
\|Ta\|_{L^{s}}\leq C \|a\|_{L^{s}}$.  Then
 $$ \|Ta\|_{L^{p}_{|x|^{\alpha}}}\leq C
  $$
 for all  $a$, the central R-$(  p,s,
\alpha)-$blocks    and    central R-$(  p,s, \alpha)-$blocks  of
restrict type,   where $C$ is independent of $a$.

\par {\bf Proof} It is obvious  since $Tf$
satisfies (3.1), (3.2), (3.4) and (3.5), see \cite{LY}.

\par Let $FRL^{p,s}_{|x|^{\alpha}}$ be the set of all finite linear
combination of R-$(p,s,\alpha)$-blocks, and
$F\dot{R}L^{p,s}_{|x|^{\alpha}}$ be the set of all finite linear
combination of R-$(p,s,\alpha)$-blocks of restrict type. We have

{\bf Corollary 3.2} Let $1\leq s < \infty,0<p\leq s,- n < \alpha <
n(p-1)  $.
 Suppose  that a sublinear operator $T$ satisfies (3.6). If  $ \|Ta\|_{L^{s}}\leq C \|a\|_{L^{s}}$,  then
 $$ \|Tf\|_{L^{p}_{|x|^{\alpha}}}\leq C \|f\|_{RL^{p,s}_{|x|^{\alpha}}}
 \eqno{(3.7)}$$
 for all $f \in FRL^{p,s}_{|x|^{\alpha}}$, and
$$ \|Tf\|_{L^{p}_{|x|^{\alpha}}}\leq C \|f\|_{\dot{R}L^{p,s}_{|x|^{\alpha}}}
 \eqno{(3.8)}$$
 for all $f \in F\dot{R}L^{p,s}_{|x|^{\alpha}}$,
 where $C$ are  independent of $f$.

\par {\bf Proof}
Let $f$ be an element of $F\dot{R}L^{p,s}_{|x|^{\alpha}}$ and pick a
representation of $f=\sum_{j=1}^N \lambda_ja_j$ such that
$$\|f\|_{\dot{R}L^{p,s}_{|x|^{\alpha}}}\approx \left(\sum\limits_{j=1}^{N}|\lambda _j|^{\bar{p}}\right)^{1/{\bar{p}}}.$$
Then (3.7) and (3.8) follow easily from Corollary 3.1.

\par {\bf Remark 3.1}  The size conditions (3.1), (3.2), (3.4),
(3.5) and (3.6) are satisfied by a lot of operators arising in
harmonic analysis. Hardy-Littlewood maximal operators, and
Littlewood-Paley $g$-function,  Lusin area function and
Littlewood-Paley $g_{\lambda}^*$-function defined as in \cite{LY}
satisfy the size conditions (3.1),  (3.2), (3.4) and  (3.5), see
\cite{LY}. While (3.6) is satisfied by many operators, such as
Calder\'{o}n-Zygmund operators, the Carleson maximal operator,
C.Fefferman's singular multiplier operator in
\cite{Hir,Wain,Fefferman70,Chanillo}, R.Fefferman's singular
integral operator in \cite{DR}, the Bochner-Riesz means at the
critical index, certain oscillatory singular integral operators in
\cite{PhongStein,RicciStein,Pan}, and so on, see \cite{LY} and
\cite{Longs}.

\par \section*
 {\bf  4 Hardy-Littlewood maximal
operators }

\par
For Hardy-Littlewood maximal operator $M$, we have more than
Theorems 3.2 and 3.3.

{\bf Theorem 4.1} Let $1< s<\infty$ and $0<p \leq s.$ Then
$$
\|Mf\|_{L^{p}_{|x|^{\alpha}}} \leq C
\|f\|_{{R}L^{p,s}_{|x|^{\alpha}}}{\rm ~~(or
~~}\|Mf\|_{L^{p}_{|x|^{\alpha}}} \leq C
\|f\|_{\dot{R}L^{p,s}_{|x|^{\alpha}}})\eqno(4.1)
$$
for all $f\in {R}L^{p,s}_{|x|^{\alpha}}$ (or $
\dot{R}L^{p,s}_{|x|^{\alpha}}$) if and only if $-n <\alpha<n(p-1 ).
$

\par {\bf Proof} For the "if" part, it is known that Hardy-Littlewood maximal
operator $M$ satisfies  (3.4) and (3.5) (see \cite{LY}).  $M$
satisfies also (3.3) for  $f=\sum\lambda_ja_j \in
{R}L^{p,s}_{|x|^{\alpha}}$ (or $\dot{R}L^{p,s}_{|x|^{\alpha}}$),
where each $a_j$ is a central R-$( p,s, \alpha)-$block (or, of
restriction type). In fact, Let $f=\sum\lambda_ja_j \in
{R}L^{p,s}_{|x|^{\alpha}}$ (or $\dot{R}L^{p,s}_{|x|^{\alpha}}$),
 taking  $x\in Q\subset {\bf R}^n$ and using
Minkowski inequality, we have
$$\frac{1}{|Q|}\int_Q|\sum\lambda_ja_j (y)|dy
\leq \sum |\lambda_j|\frac{1}{|Q|}\int_Q|a_j(y)|dy \leq \sum
|\lambda_j|Ma_j(x)$$
 for $x\in Q$, and it follows that $Mf(x) \leq \sum |\lambda_j|Ma_j(x)$. So, the "if"
 part of Theorem 4.1 is a corollary of Theorem 3.3

  \par  For
the "only if " part, we need only to prove that (4.1) fails for
$\alpha\leq -n$ and $\alpha \geq n(p-1 )$.  Let us first see the
case of  $\alpha\leq -n$. We know that
$$
M(\chi_{C_0})(x)\geq C\eqno(4.2)
$$
for $x\in B_0$, where $C_0=B_0\setminus B_{-1}$. It is clear that
$\chi_{C_0}(x) \in \dot{R}L^{p,s}_{|x|^{\alpha}}$ ( and
$RL^{p,s}_{|x|^{\alpha}}$)
 for any $0<p,s<\infty$ and $
-\infty<\alpha<\infty$ with
$\|\chi_{C_0}\|_{\dot{R}L^{p,s}_{|x|^{\alpha}}}\leq 1$ (and
$\|\chi_{C_0}\|_{RL^{p,s}_{|x|^{\alpha}}}\leq 1$ ).
\par By (4.2), we have
\begin{eqnarray*}
\int_{{\bf R}^n} | M(\chi_{C_0})(x)|^p|x|^{\alpha}dx& \geq & C
\int_{B_0}  |x|^{\alpha} dx
\\ &=& C \int_0^1  r^{\alpha+n-1}dr
\\&=& \infty
\end{eqnarray*}
for $\alpha\leq -n$.  Thus, (4.1) fails for $\alpha\leq - n $. Next,
let us consider the case of
 $\alpha\geq n(p-1)$, it suffices to prove the following proposition.

\par {\bf Proposition 4.1} Let $0<p<\infty, n(p-1)\leq \alpha <\infty,$ and
$0<q\leq 1, f\in L^q$. Then $Mf$ is never in ${L}^{p}
  _{|x|^{\alpha}}$ if $ f  \neq 0.$

\par Once Proposition 4.1 is proved, by Proposition 2.2i, then (4.1) fails for $\alpha\geq
n(p-1)$ (i.e. $\frac{np}{n+\alpha}\leq 1$).  Thus,   Theorem 4.1
holds.

\par{\bf Proof of Proposition 4.1} If supp $f \subset B(0,R) $, then
$$Mf(x)\geq \frac{\|f\|_{L^1}}{v_n}\frac{1}{(|x|+R)^n},~~{\rm for }~~ |x|\geq R, $$
where $v_n$ is the volume of the unit ball (see \cite{Grafakos}). By
H\"{o}lder inequality, it follows that
$$Mf(x)\geq \frac{\|f\|_{L^q}}{v_nR^{(1-q)/q}}\frac{1}{(|x|+R)^n},~~{\rm for }~~ |x|\geq R. $$
Noticing $np-\alpha\leq n$, we have
\begin{eqnarray*}
\int |Mf(x)|^p|x|^\alpha dx &\geq& \left(
\frac{\|f\|_{L^q}}{v_nR^{(1-q)/q}}\right)^{p}\int _{|x|\geq
R}\frac{|x|^\alpha}{(|x|+R)^{np}}dx
\\&
\simeq & \left( \frac{\|f\|_{L^q}}{v_nR^{(1-q)/q}}\right)^{p}\int
_{|x|\geq R}\frac{1}{(|x|+R)^{np-\alpha}}dx,
\\&
\geq & \left( \frac{\|f\|_{L^q}}{v_nR^{(1-q)/q}}\right)^{p}\int
_{|x|\geq R}\frac{1}{(|x|+R)^{n}}dx=\infty,
\end{eqnarray*}
if $\|f\|_{L^q}\neq 0.$ So if $Mf$ is in $ {L}^{p}
  _{|x|^{\alpha}}$, then $f=0$ a.e.. For general $f$, taking $f_R(x)=f(x)\chi_{\{|x|\leq
  R\}}$, then $f_R(x)=0 $ for almost all $x$ in the ball of radius
  $R>0.$ Thus, $f=0 $ a.e.. We finish the proof of Proposition 4.2.

\par {\bf Theorem 4.2}  Let $s=1, 0<p\leq 1, -\infty <\alpha \leq n(p-1).$
Then there exists $f\in \dot{R}L^{p,s}_{|x|^\alpha}$ such that (4.1)
fails.

\par {\bf Proof} When  $s=1, 0<p\leq 1$ and $ -\infty <\alpha \leq n(p-1),$
 by Theorem 2.1 ii), we
have $L^1\subset \dot{R}L^{p,s}_{|x|^\alpha}$. Then, by Proposition
4.2, (4.1) fails for $f\in L^1.$

\par
\section*
{\bf 5 Singular integrals }

\par Let us first consider Calder\'{o}n-Zygmund
operator $T$. Let $ T_\varepsilon $ and $T^* $ denote the
corresponding truncated operators and maximal operators
respectively. See \cite{Mey} for the definitions. We know that  $T,
T_\varepsilon $ and $T^* $ are bounded on  $(L^p)$ for $1<p<\infty$.
 From the good-$\lambda$ inequality, we have

\par {\bf Lemma 5.1} Let $0<p<\infty, \omega \in A_\infty, $ and $T^*$ be a maximal Calder\'{o}n-Zygmund
operator. Then
 $$\|T^*f\|_{L^{p} _{|x|^{\alpha}}}\leq
C \|Mf\|_{L^{p} _{|x|^{\alpha}}}
 $$
 for all $f\in L^1_{loc}$, where
$C$ are independent of $f$.

\par See \cite{Mey} for a proof of Lemma 5.1.
\par {\bf Theorem 5.1}   Let $1<s<\infty,  0< p \leq s ,-n < \alpha < n(p-1).
$ Then,  for the functions in $RL^{p,s} _{|x|^{\alpha}}$ and
$\dot{R}L^{p,s} _{|x|^{\alpha}}$,
 $Tf, T_\varepsilon f$ and $T^* f$  are  defined for
 $\mu_\alpha$-almost every points,
 $T, T_\varepsilon $ and $T_* $  are bounded from $RL^{p,s}
_{|x|^{\alpha}}$ into $L^{p} _{|x|^{\alpha}}$ and from
$\dot{R}L^{p,s} _{|x|^{\alpha}}$ into $L^{p} _{|x|^{\alpha}}$, and
$$
T_\varepsilon f=\sum _i \lambda _i T_\varepsilon a_i, ~~\mu_\alpha
{\rm-
  a.e. } \eqno(5.1)
$$
$$
T f=\sum _i \lambda _i T a_i, ~~\mu_\alpha {\rm-
  a.e. } \eqno(5.2)
$$
and
$$T^*f (x)\leq \sum _i |\lambda
_i|T^* a_i(x),~\mu_\alpha {\rm-
  a.e. }\eqno(5.3)$$
for
 $f=\sum _i \lambda _i a_i \in RL^{p,s}_{|x|^\alpha}
 $ and $\dot{R}L^{p,s}_{|x|^\alpha}  $,
 where $a_i$ are central R-$(p,s,\alpha) $ - blocks
 (or central R-$(p,s,\alpha) $ - blocks of restrict type), and $\sum_i
 |\lambda_i|^{\bar{p}} <+\infty$.

\par {\bf Proof} From Theorem 4.1 and Lemma 5.1, we get that
  $T, T_\varepsilon $ and $T^* $ are bounded   from $RL^{p,s}
_{|x|^{\alpha}}$ into $L^{p} _{|x|^{\alpha}}$ and from
$\dot{R}L^{p,s} _{|x|^{\alpha}}$ into $L^{p} _{|x|^{\alpha}}$. It
follows that $Tf, T_\varepsilon f$ and $T^* f$  are  defined for the
functions $f $ in $RL^{p,s} _{|x|^{\alpha}}$ and $\dot{R}L^{p,s}
_{|x|^{\alpha}}$ for
 $\mu_\alpha$-almost every points.

Next, let us  prove  (5.1). Firstly we consider the case of
$RL^{p,s} _{|x|^{\alpha}}$.
 Let  $f=\sum _{i=0}^\infty \lambda _i a_i
 $ where $a_i$ are central R-$(p,s,\alpha) $ - block and $\sum_{i=0}^\infty
 |\lambda_i|^{\bar{p}} <+\infty.$ For any $\varepsilon>0$, there
 exits an  $N_0>0$ such that
 $$\sum_{i=N_0+1}^\infty
 |\lambda_i|^{\bar{p}} <\varepsilon.\eqno(5.4)$$
By  the sublinearity of the quantity
$\|\cdot\|_{RL^{p,s}_{|x|^{\alpha}}}$, and noticing that
$\|a_i\|_{RL^{p,s}_{|x|^\alpha}}<1$, we have
$$\|f-\sum _{i=0}^{N_0} \lambda _i a_i\|_{RL^{p,s}_{|x|^\alpha}}
=\|\sum _{i=N_0+1}^{\infty} \lambda _i
a_i\|_{RL^{p,s}_{|x|^\alpha}}\leq  \sum _{i=N_0+1}^{\infty} |\lambda
_i|\leq  \left\{\sum _{i=N_0+1}^{\infty} |\lambda
_i|^{\bar{p}}\right\}^{1/\bar{p}}<\varepsilon^{1/\bar{p}}.\eqno(5.5)$$
 Then, by the boundedness result of $T_\varepsilon$ from $RL^{p,s}
_{|x|^{\alpha}}$ into $L^{p} _{|x|^{\alpha}}$, (5.5), (5.4) and $
\|T_\varepsilon a_i \|_{L^{p }_{|x|^\alpha}}\leq C
\|a_i\|_{RL^{p,s}_{|x|^\alpha}}\leq C  $, we have
\begin{eqnarray*}
\|T_\varepsilon f-\sum _{i=0}^{\infty} \lambda _iT_\varepsilon
a_i\|^{\bar{p}}_{L^{p }_{|x|^\alpha}}
& \leq & \|T_\varepsilon
(f-\sum _{i=0}^{N_0} \lambda _i  a_i)\|^{\bar{p}}_{L^{p
}_{|x|^\alpha}} +\|\sum _{i=N_0+1}^{\infty} \lambda _i T_\varepsilon
a_i\|^{\bar{p}}_{L^{p }_{|x|^\alpha}}
\\ & \leq & \|
f-\sum _{i=0}^{N_0} \lambda _i
a_i\|^{\bar{p}}_{RL^{p,s}_{|x|^\alpha}} +\sum _{i=N_0+1}^{\infty}
|\lambda _i|^{\bar{p}} \|T_\varepsilon a_i\|^{\bar{p}}_{L^{p
}_{|x|^\alpha}}
\\&<& C  \varepsilon .
\end{eqnarray*}
Let $\varepsilon \rightarrow 0$, we get that $ \|T_\varepsilon
f-\sum _{i=0}^{\infty} \lambda _iT_\varepsilon a_i
\|_{{L}^{p}_{|x|^\alpha}} =0,
 $ and it follows that $
T_\varepsilon f=\sum _{i=0}^{\infty} \lambda _iT_\varepsilon a_i, $
$ ~\mu_\alpha {\rm -
  a.e. }.
 $ Thus, we have proved (5.1).

  \par   (5.2) can be proved similarly. (5.3) follows from (5.1) easily. The case of $\dot{R}L^{p,s}
_{|x|^{\alpha}}$ is similar. Thus, we finish the proof of Theorem
5.1.

It is known that Hilbert transform, Riesz transform   and the
regular singular integral operators   defined in \cite{GRubio}
  are Calder\'{o}n-Zygmund operators. So

\par {\bf Corollary 5.1} For Hilbert transform, Riesz transform, the
regular singular integral operators, and their corresponding
truncated operators and maximal operators, the same conclusions hold
as those stated in Theorem 5.1.

\par When the index exceed the critical index, Bochner-Riesz means
are controlled pointwise by Hardy-Littlewood maximal operator. So it
is easy to see  from the theorem above that

\par {\bf Corollary 5.2} For Bochner-Riesz means with the index exceeding the critical index,
the same conclusions hold
as those stated in Theorem 5.1.

\par For Hilbert transform we have more.

\par {\bf Theorem 5.2} Let $1< s<\infty$ and $0<p \leq s.$ Then $H$   is
bounded from $RL^{p,s}_{|x|^{\alpha}}({\bf R})$ or
$\dot{R}L^{p,s}_{|x|^{\alpha}}({\bf R})$ into $
L^{p}_{|x|^{\alpha}}({\bf R})$, i.e.
$$
\|Hf\|_{L^{p}_{|x|^{\alpha}}({\bf R})} \leq C
\|f\|_{RL^{p,s}_{|x|^{\alpha}}({\bf R})},\eqno(5.6)
$$
or
$$
\|Hf\|_{L^{p}_{|x|^{\alpha}}({\bf R})} \leq C
\|f\|_{\dot{R}L^{p,s}_{|x|^{\alpha}}({\bf R})},\eqno(5.7)
$$
 if and only if $-1  <\alpha<p-1 . $  So is $H^*$.

\par {\bf Proof} By Corollary 5.1, we need only to prove that (5.6) and (5.7) fail for $\alpha\leq -1$
and $\alpha \geq p-1 $. Let $f=\chi_{[1,2]}$, it is easy to see
 that $f\in RL^{p,s}_{|x|^{\alpha}}({\bf R})$ and  $  \dot{R}L^{p,s}_{|x|^{\alpha}}({\bf R})$, and
$$
Hf(x)=\frac{1}{\pi}\log \left|\frac{x-1}{x-2}\right|.
$$
It is clear that $Hf(x)\sim \frac{1}{\pi (x-2)}$ when
$|x|\rightarrow \infty$ or $|x|\rightarrow 0$. It follows that $Hf$
is  not in $L^{p}_{|x|^{\alpha}}({\bf R})$ as   $\alpha \geq p-1 $,
and  clearly, neither is $H^*$. Since $Hf(x)$ is continues and not
equal to zero on $[0,1/2]$, it is easy to see that $Hf$ is not in
$L^{p}_{|x|^{\alpha}}({\bf R})$ as $\alpha\leq -1$. $H^*$ follows.
Thus, we finish the proof of Theorem 5.2.

\par
Next, we will also see that many other linear operators  defined for
the functions in $L^p$ with $1<p<\infty$ have  extensions which map
$ RL^{p,s}_{|x|^\alpha}$ ( or $ \dot{R}L^{p,s}_{|x|^\alpha}$) into $
{L}^{p}_{|x|^\alpha}$.

{\bf Theorem 5.3} Let $1<s<\infty,  0< p \leq s ,-n < \alpha <
n(p-1)$.
  Suppose that    $T$ is a bounded linear operator on  $L^s$, which satisfies (3.6).
 Then
 \par i)
  $T$ has an extension which maps     $
RL^{p,s}_{|x|^\alpha}$ into $ {L}^{p}_{|x|^\alpha}$ and  satisfies
(5.2) and
$$\|T f\|_{L^{p} _{|x|^{\alpha}}}\leq
C \|f\|_{ RL^{p,s} _{|x|^{\alpha}}} \eqno(5.8)
 $$
 for  $f=\sum\lambda_ja_j \in RL^{p,s} _{|x|^{\alpha}}$, where
each $a_j$ is a central R-$(  p,s, \alpha)-$block.

\par ii) The conclusions above hold for
        $ \dot{R}L^{p,s}_{|x|^\alpha}$.

\par {\bf Proof} Let us first consider the case of
$RL^{p,s} _{|x|^{\alpha}}$.
 Let $f\in RL^{p,s}_{|x|^\alpha}$, then $f=\sum _i \lambda _i a_i
 $ where $a_i$ are central R-$(p,s,\alpha) $ - block and $\sum_i
 |\lambda_i|^{\bar{p}} <+\infty.$ By Corollary 3.1, we have
  $$
\|\sum _i \lambda _i T a_i\|^{\bar{p}}_{L^{p }_{|x|^\alpha}}
 \leq   \sum _i |\lambda _i |^{\bar{p}} \|T  a_i\|^{\bar{p}}_{L^{p }_{|x|^\alpha}}
  \leq C \sum _i |\lambda _i |^{\bar{p}} <\infty.
  $$
It follows that $|\sum _i \lambda _i T  a_i(x)|<\infty, ~\mu_\alpha
{\rm-
  a.e. }$. Set
  $$T  f(x) =\sum _i \lambda _i T  a_i(x) ,
~\mu_\alpha {\rm-
  a.e. }.\eqno (5.9)$$
 Once it is proved that (5.9) holds for all  decomposition   of $f$
  as $f=\sum _{i=0}^\infty \lambda _i a_i
 $ where $a_i$ are central R-$(p,s,\alpha) $ - block and $\sum_{i=0}^\infty
 |\lambda_i|^{\bar{p}} <+\infty,$ then the extension required is
 obtained.

  Let
$$
f=\sum _i \lambda _i^{(1)} a_i^{(1)}=\sum _i \lambda _i^{(2)}
a_i^{(2)}  \eqno (5.10)
 $$
 with
$$\sum_i
 |\lambda_i^{(1)}|^{\bar{p}} <+\infty
 {~~\rm and ~~} \sum_i
 |\lambda_i^{(2)}|^{\bar{p}} <+\infty \eqno (5.11)$$
 and  $a_i^{(1)}$ and  $a_i^{(2)}$ are central R-$(p,s,\alpha) $ - blocks.
Once $  \sum _i \lambda _i^{(1)} Ta_i^{(1)}=\sum _i \lambda _i^{(2)}
Ta_i^{(2)},~\mu_\alpha {\rm -
  a.e. }
 $, then we get the extension required. In fact, for any $\delta>0$, by (5.11), there exists an $i_0$ such that
$$\sum_{i=i_0}^\infty
 |\lambda_i^{(1)}|^{\bar{p}} <\delta^{\bar{p}}
 {~~\rm and ~~} \sum_{i=i_0}^\infty
 |\lambda_i^{(2)}|^{\bar{p}} <\delta^{\bar{p}}. \eqno (5.12)$$
From (5.10), we see that
$$
\sum _{i=1}^{i_0-1} (\lambda _i^{(1)} a_i^{(1)}- \lambda _i^{(2)}
a_i^{(2)} )= \sum _{i=i_0}^{\infty}\lambda _i^{(2)} a_i^{(2)}-\sum
_{i=i_0}^{\infty} \lambda _i^{(1)} a_i^{(1)},
 $$
then,
$$
\|\sum _{i=1}^{i_0-1} (\lambda _i^{(1)} a_i^{(1)}- \lambda _i^{(2)}
a_i^{(2)} )\|^{\bar{p}}_{RL^{p,s}_{|x|^\alpha}} \leq
\sum_{i=i_0}^\infty
 |\lambda_i^{(1)}|^{\bar{p}} + \sum_{i=i_0}^\infty
 |\lambda_i^{(2)}|^{\bar{p}} <2 \delta^{\bar{p}}. \eqno (5.13)
 $$
By the linearity of $T $, we have
$$
\sum _{i=1}^{\infty} \lambda _i^{(1)} T  a_i^{(1)}- \sum
_{i=1}^{\infty} \lambda _i^{(2)} T  a_i^{(2)}
 =T  (\sum
_{i=1}^{i_0-1} (\lambda _i^{(1)} a_i^{(1)}- \lambda _i^{(2)}
a_i^{(2)} ))+ \sum _{i=i_0}^{\infty}\lambda _i^{(1)} T
a_i^{(1)}-\sum _{i=i_0}^{\infty} \lambda _i^{(2)} T  a_i^{(2)}.\eqno
(5.14)
 $$
By Corollary 3.2 and (5.13), we see that
$$
\|T  (\sum _{i=1}^{i_0-1} (\lambda _i^{(1)} a_i^{(1)}- \lambda
_i^{(2)} a_i^{(2)} ))\|^{\bar{p}}_{{L}^{p}_{|x|^\alpha}} \leq \|\sum
_{i=1}^{i_0-1} (\lambda _i^{(1)} a_i^{(1)}- \lambda _i^{(2)}
a_i^{(2)} )\|^{\bar{p}}_{RL^{p,s}_{|x|^\alpha}}<2 \delta^{\bar{p}}.
\eqno (5.15)
 $$
From (5.14), (5.15), Corollary 3.1 and  (5.12), we have
$$
\|\sum _{i=1}^{\infty} \lambda _i^{(1)} T  a_i^{(1)}- \sum
_{i=1}^{\infty} \lambda _i^{(2)} T  a_i^{(2)}
\|^{\bar{p}}_{{L}^{p}_{|x|^\alpha}} <4 \delta^{\bar{p}}.
 $$
Let $\delta\rightarrow 0$, we get that $ \|\sum _{i=1}^{\infty}
\lambda _i^{(1)} T  a_i^{(1)}- \sum _{i=1}^{\infty} \lambda _i^{(2)}
T a_i^{(2)} \|^{\bar{p}}_{{L}^{p}_{|x|^\alpha}} =0.
 $ It follows that $
\sum _{i=1}^{\infty} \lambda _i^{(1)} T  a_i^{(1)}=\sum
_{i=1}^{\infty} \lambda _i^{(2)} T  a_i^{(2)}, ~\mu_\alpha {\rm -
  a.e. }
 $. Thus, (5.9) holds for all  decomposition   of $f$
  as $f=\sum _{i=0}^\infty \lambda _i a_i \in RL^{p,s}
  _{|x|^{\alpha}},
 $ i.e. (5.2) holds for the extension of $T$.


  (5.8) follows from (5.2) easily. The proof of the case of  $
  \dot{R}L^{p,s}_{|x|^\alpha}$ is similar.
Thus, we finish the proof of Theorem 5.3.

\par {\bf Corollary 5.3} For  C.Fefferman's singular multiplier
operator, R.Fefferman's singular integral operator,  Bochner-Riesz
means at the critical index and certain oscillatory singular
integral operators, the same conclusions hold as those stated in
Theorem 5.3.

\section*
{\bf  6 Fourier series    }

\par For $N>0$, define the Dirichlet summation operator
$$
S_N f (x)= \int_{-N}^{N} \hat{f}(\xi)e^{2 \pi i \xi \cdot x}d\xi
\eqno(6.1)
$$
on $\mathcal{S}({\bf R})$. We have the identity
$$S_Nf=\frac{i}{2}(  M^NH M^{-N}f-M^{-N}H  M^{N}f) \eqno(6.2)$$
for $f\in\mathcal{S}({\bf R})$, where $H$ is Hilbert transform and
$M^{N}g(x)=e^{2\pi i Nx}g(x)$. By (6.2), $S_N$ extends to a bounded
operator on $L^p({\bf R})$ with $1<p<\infty$.  By Theorem 5.3,   we
have
\par {\bf Theorem 6.1}  Let $  1<s<\infty,  0<p\leq s,-1<\alpha<p-1 .$
Then, the operator $S_N$, initially defined for $f\in
\mathcal{S}({\bf R})$  by (6.1),  extends to a bounded operator from
$RL^{p,s}_{|x|^\alpha}$ (or $\dot{R}L^{p,s}_{|x|^\alpha}$) into
$L^p_{|x|^\alpha}({\bf R})$, and
$$ S_Nf=\sum _i \lambda _i S_N a_i,
~~\mu_\alpha {\rm-
  a.e. }\eqno(6.3)
$$
 and
$$
\|S_N f\|_{L^{p }_{|x|^\alpha}({\bf R})}\leq C
\|f\|_{RL^{p,s}_{|x|^\alpha}({\bf R})}, ({\rm ~or ~}\|S_N f\|_{L^{p
}_{|x|^\alpha}({\bf R})}\leq C
\|f\|_{\dot{R}L^{p,s}_{|x|^\alpha}({\bf R})}),\eqno(6.4)
$$
for
 $f=\sum _i \lambda _i a_i \in RL^{p,s}_{|x|^\alpha}
 $ (or $\dot{R}L^{p,s}_{|x|^\alpha}$),  where $a_i$ are central R-$(p,s,\alpha) $ - blocks
 (or central R-$(p,s,\alpha) $ - blocks of restrict type) and $\sum_i
 |\lambda_i|^{\bar{p}} <+\infty.$

 \par
 Define  Carleson operator
$$
C f (x)= \sup _{N>0}\left|S_Nf(x)\right|\eqno(6.5)
$$
for $RL^{p,s}_{|x|^\alpha}({\bf R})$ (or
$\dot{R}L^{p,s}_{|x|^\alpha}({\bf R})$). It is easy to see that the
conclusion of  Corollary 3.1 holds for Carleson operator $C$. We
have

\par {\bf Theorem 6.2}  Let $  1<s<\infty,  0<p\leq s,-1<\alpha<p-1 .$
Then,
$$ Cf\leq \sum _i \lambda _iC a_i,
~~\mu_\alpha {\rm-
  a.e. }\eqno(6.6)
$$
 and
$$
\|C f\|_{L^{p }_{|x|^\alpha}({\bf R})}\leq C
\|f\|_{RL^{p,s}_{|x|^\alpha}({\bf R})} ({\rm ~or ~}\|C f\|_{L^{p
}_{|x|^\alpha}({\bf R})}\leq C
\|f\|_{\dot{R}L^{p,s}_{|x|^\alpha}({\bf R})}),\eqno(6.7)
$$
 for
 $f=\sum _i \lambda _i a_i \in RL^{p,s}_{|x|^\alpha}
 $ (or $\dot{R}L^{p,s}_{|x|^\alpha}$), where $a_i$ are central R-$(p,s,\alpha) $ - blocks
 (or central R-$(p,s,\alpha) $ - blocks of restrict type), and $\sum_i
 |\lambda_i|^{\bar{p}} <+\infty.$
\par {\bf Proof} (6.6) follows from (6.3). (6.7) follows from (6.6)
and Corollary 3.1.

\par According to a well-known patten, we can get  the $\mu_\alpha$-a.e. convergence
 and $L^p({\bf R}^n)$-norm convergence of $S_Nf$ for $f\in
RL^{p,s}_{|x|^\alpha}$ and $ RL^{p,s}_{|x|^\alpha}$ from Theorems
6.1 and 6.2 respectively. In fact, the convergence results follow
from the following Lemmas whose proofs are similar to the cases of
$L^p$.

\par {\bf Lemma 6.1 (Uniform boundedness principle)}
Let $1<s<\infty, 0<p\leq s, -n<\alpha<n(p-1)$.

\par (i) Let $\mathcal{D}$ be a dense subspace of $
RL^{p,s}_{|x|^\alpha}$ and suppose that $T_R$ is a sequence of
linear operators such that
$$
\|T_Rf-f\|_{L^p_{|x|^\alpha}} \rightarrow 0, R\rightarrow \infty,
$$
 for $f\in \mathcal{D}$.
 Then in order for $T_Rf\rightarrow f$ in
$L^p_{|x|^\alpha}$ norm  as $R\rightarrow \infty $ for all $f\in
RL^{p,s}_{|x|^\alpha}$ it is a necessary and sufficient condition
that we have the estimate
$$
 \|T_Rf\|_{L^p_{|x|^\alpha}}\leq C\|f\|_{RL^{p,s}_{|x|^\alpha}}
$$
for all sufficiently large $R,$ where the constants $C$ are
independent of $R$.

\par (ii)  For $\dot{R}L^{p,s}_{|x|^\alpha}$,  the same
conclusions hold as those stated in (i).

\par Let $T_\varepsilon$ be a linear operator for every
$\varepsilon>0$, and
$T_*(f)(x)=\sup_{\varepsilon>0}|T_\varepsilon(f)(x)|.$
\par {\bf Lemma 6.2 (Maximal principle)} Let $0<p<\infty,1\leq s<\infty,
-\infty<\alpha<\infty$.

\par (i) Let $\mathcal{D}$ be a dense subspace of
$RL^{p,s}_{|x|^\alpha}$, and suppose that for some $C$ and for all
$f\in RL^{p,s}_{|x|^\alpha}$ we have
$$
 \|T_*f\|_{L^p_{|x|^\alpha}}\leq C\|f\|_{RL^{p,s}_{|x|^\alpha}}
 $$
and for all $f\in \mathcal{D}$
$$
 \lim_{\varepsilon \rightarrow 0}T_\varepsilon(f)(x)=f(x) ~~~~~\mu_\alpha
 {\rm-
 a.e..}\eqno(6.8)
 $$
Then for all $f\in RL^{p,s}_{|x|^\alpha}$ (6.8) holds.
\par (ii)  For $\dot{R}L^{p,s}_{|x|^\alpha}$,  the same
conclusions hold as those stated in (i).


\par It is easy to see that a.e. convergence implies
$\mu_\alpha$-a.e. convergence when $\mu_\alpha \in A_\infty$. Then $
\lim_{N\rightarrow \infty}S_N f=f,\mu_\alpha {\rm-
  a.e.}$,
for all $ f\in L^s({\bf R})$ with $  1<s<\infty$ and $
-1<\alpha<\infty,$ since $ \lim_{N\rightarrow \infty}S_N f=f  ~{\rm
  a.e.}$ for all $ f\in L^s({\bf R})$ with $1<s<\infty$.

\par On the other hand, we see that $|S_N f|\leq C f $. By Theorem 6.2, it follows that
 $|S_N f-f|\leq C f +|f| \in L^p_{|x|^\alpha}({\bf R})$
 for  $ f\in L^s({\bf R})\cap RL^{p,s}_{|x|^\alpha}({\bf
R})$ (or $  L^s({\bf R})\cap \dot{R}L^{p,s}_{|x|^\alpha}({\bf R})$)
if $1<s<\infty,  0<p\leq s$ and $ -1<\alpha< p-1 $. The Lebesgue
dominated
  convergence theorem gives that $S_N f$ converge to $f$ in
  $L^p_{|x|^\alpha}({\bf R})$.

\par It is easy to see that $  L^s({\bf R})\cap RL^{p,s}_{|x|^\alpha}({\bf
R})$ in $   RL^{p,s}_{|x|^\alpha}({\bf R})$ and $  L^s({\bf R})\cap
\dot{R}L^{p,s}_{|x|^\alpha}({\bf R})$ in $
\dot{R}L^{p,s}_{|x|^\alpha}({\bf R})$ are dense respectively. Thus,
by the uniform boundedness principle   and the maximal principle
above, we obtain that

\par {\bf Theorem 6.3} Let   $ 1< s<\infty, 0<p\leq s $ and $-1<\alpha<p-1.$
Then
$$ \lim_{N\rightarrow \infty}S_N f=f~~~~{\rm
  in } ~~L^p_{|x|^\alpha}({\bf R})$$
and
$$ \lim_{N\rightarrow \infty}S_N f=f,~~~~\mu_\alpha {\rm-
  a.e.}$$
  for all
$f\in RL^{p,s}_{|x|^\alpha}({\bf R})\cup
\dot{R}L^{p,s}_{|x|^\alpha}({\bf R})$ .

\par For the the functions on $I=[-1,1]$, Theorem 6.3 implies some
 weighted result. Let us first prove the following
proposition.

\par {\bf Proposition 6.1} (i) Let $1\leq  s< \infty, 0<p< s,
-n<\alpha <\infty$.  Then
$$
{L}^{s}_{|x|^{\alpha}}(B_0)\subset RL^{p,s}_{|x|^{\alpha}}(B_0).
$$
(ii) Let $1\leq  s< \infty, 0<p<\infty, -n<\alpha <\infty$,  then
$$
{L}^{s}_{|x|^{\alpha}}(B_0)\subset
\dot{R}L^{p,s}_{|x|^{\alpha}}(B_0).
$$
Here, $B_0=\{x\in {\bf R}^n:|x|\leq 1\}$.
\par {\bf Proof } (i)
  Let $f\in L^{s}_{|x|^\alpha}(B_0),$   write
$f=\sum_{k=-\infty }^{0} \lambda_k b_k,$ where
$$
b_k=\frac{f\chi_{C_k}}{|B_k|^{\frac{\alpha}{pn}+\frac{1}{p}-\frac{1}{s}}\|f\chi_{C_k}\|_{L^{s}(B_0)}},$$
$$
\lambda_k=
|B_k|^{\frac{\alpha}{pn}+\frac{1}{p}-\frac{1}{s}}\|f\chi_{C_k}\|_{L^{s}(B_0)}\approx
|B_k|^{ (\frac{1}{p}-\frac{1}{s})
(1+\frac{\alpha}{n})}\|f\chi_{C_k}\|_{L^{s}_{|x|^{\alpha} }(B_0)}
$$
$k=0,-1,-2,\cdots,$  $B_k=B(0,2^k), C_k=B_k\setminus B_{k-1}. $ We
see that all $b_k$ are $(p,s,\alpha)-$ blocks, and by H\"{o}lder
inequality as $\bar{p}<s$,

\begin{eqnarray*}
 \sum_{k=-\infty }^{0} |\lambda_k|^{\bar{p}} &\leq&  C \sum_{k=-\infty }^{0} 2^{kn\bar{p} (\frac{1}{p}-\frac{1}{s})
(1+\frac{\alpha}{n})}\|f\chi_{C_k}\|^{\bar{p}}
_{L^{s}_{|x|^{\alpha}(B_0)}}
\\ &\leq& C
\left(\sum_{k=-\infty }^{0}2^{kn (\frac{1}{p}-\frac{1}{s})
(1+\frac{\alpha}{n})\frac{s}{s-\bar{p}}}\right)^{(s-\bar{p})/s}\left(\sum_{k=-\infty
}^{0}
\|f\chi_{C_k}\|^s_{L^{s}_{|x|^{\alpha}}(B_0)}\right)^{\bar{p}/s}
\\&= &
C\|f\|^{\bar{p}}_{L^{s}_{|x|^{\alpha}}(B_0)},
\end{eqnarray*}
noticing that $p<s$ and $ -n <\alpha.$ Thus, (i) holds.
\par (ii) Each $f\in L^{s}_{|x|^\alpha}(B_0)$ is a  central R-$ (p, s, \alpha)$-block of restrict
type, then (ii) is follows.

\par From Theorem 6.3 and Proposition 6.1 (ii), it is easy to see that

\par {\bf Corollary 6.1} Let   $I=[-1,1], 1< s<\infty  $ and $ -1<\alpha<
s-1.$ Then
$$ \lim_{N\rightarrow \infty}S_N f=f~~~~{\rm
  in } ~~L^s_{|x|^\alpha}(I)$$
and
$$ \lim_{N\rightarrow \infty}S_N f=f,~~~~\mu_\alpha {\rm-
  a.e.}$$
  for all
 $f \in L^s_{|x|^\alpha}(I)$.

\par See also \cite{HY} for Corollary 6.1.

\par It  should be pointed out that
Theorem 6.3 does not imply new information  for the pointwise
convergence when $ p/s-1<\alpha< p-1, $ since
$RL^{p,s}_{|x|^{\alpha}}\subset L^{\frac{p}{1+\alpha}}
 $ and $\frac{p}{1+\alpha}>1$ by Proposition 2.2i.

\par {\bf Theorem 6.4}
 If
 $  s=1, 0<p\leq 1, -1 <\alpha \leq p-1. $
  Then there exists $f\in RL^{p,s}_{|x|^\alpha}({\bf R})$ such that
$$
\limsup_{N\rightarrow \infty}S_Nf(x)=\infty,~~\mu_\alpha{\rm -a.e.}.
\eqno(6.9)
$$

\par {\bf Proof}
 When $  s=1, 0<p\leq 1, -1 <\alpha \leq p-1, $ by Proposition 2.2ii,
we see that $L^1\subset RL^{p,s} _{|x|^\alpha}.$ Then Kolmogorov's
example shows that there exists $f\in  RL^{p,s} _{|x|^\alpha}$ such
that
$$
\limsup_{N\rightarrow \infty}S_Nf(x)=\infty~~{\rm a.e.},
$$
then (6.9) follows.

\par
Shunchao Long
\par Department of Mathematics,
\par Xiangtan University,
 \par Hunan, 411105  P.R.China
 \par E-mail: sclong@xtu.edu.cn

\end{document}